\documentclass{amsart}
\usepackage{amscd,amsthm}
\usepackage{latexsym}
\usepackage{capt-of}
\usepackage{graphicx}

\DeclareFontFamily{U}{wncy}{}
\DeclareFontShape{U}{wncy}{m}{n}{<->wncyr10}{}
\DeclareSymbolFont{mcy}{U}{wncy}{m}{n}
\DeclareMathSymbol{\Sha}{\mathord}{mcy}{"58}

\usepackage{amssymb,amsmath}
\usepackage{amsfonts}
\usepackage{comment}
\usepackage{hyperref}
\newcounter{ctfig}

\newcommand{\C}{\mathcal{C}}

\newcommand{\F}{\mathcal{F}}

\newcommand{\JacC}{{\hbox{Jac}_{\lower.5pt\hbox{$_\C$}}}}
\newcommand{\JacF}{{\hbox{Jac}_{\lower.5pt\hbox{$_\F$}}}}

\theoremstyle{plain}
\newtheorem{thm}{Theorem}[section]
\newtheorem{conj}[thm]{Conjecture}

\theoremstyle{definition}

\def\F{{\mathbb F}}

\def\C{{\mathbb C}}

\def\e32{{{}_3E_2}}
\def\f32{{{}_3F_2}}
\def\a32{{{}_3A_2}}
\makeatletter
\@namedef{subjclassname@2020}{%
  \textup{2020} Mathematics Subject Classification}
\makeatother

\begin{document}
\bibliographystyle{plain}
\bibstyle{plain}

\title[ECAAC]{A question of Erd\"{o}s on $3$-powerful numbers and an elliptic curve analogue of the Ankeny-Artin-Chowla conjecture.}

\author{P.G. Walsh}
\address{Department of Mathematics\\
University of Ottawa}
\email{gwalsh@uottawa.ca}

\date{\today}

\begin{abstract}
We describe how the Mordell-Weil group of rational points on a certain family of elliptic curves give rise to solutions to a conjecture of Erd\"{o}s on $3$-powerful numbers, and state a related conjecture which can be viewed as an elliptic curve analogue of the famous Ankeny-Artin-Chowla conjecture.
\end{abstract}

\maketitle

\section{Introduction}
In a series of papers \cite{Er1},\cite{Er2},\cite{Er3}, Erd\"{o}s posed a number of problems concerning powerful numbers. Quite recently, the arithmetic of elliptic curves has been used by Bajpai, Bennett and Chan \cite{BBC} to prove the existence of infinitely many quadruples of four pairwise coprime powerful numbers in arithmetic progression, thereby answering one of Erd\"{o}s' problems. In this article we consider a different problem of Erd\"{o}s. In particular, he asked for solutions to $a+b=c$ in coprime $3$-powerful numbers. This was solved by Nitaj \cite{Ni}, and later in a different way by Cohn \cite{Cohn}. We will describe how one can produce solutions to Erd\"{o}s' problem by using the group of rational points on an elliptic curve. It is worth noting here that despite the abundance of coprime solutions to $a+b=c$ in $3$-powerful numbers, the {\em abc} conjecture implies that there are only finitely many triples $(a,b,c)$ consisting of pairwise coprime $4$-powerful numbers which satisfy $a+b=c$, and that no such example is even known to exist (see \cite{KL} for details).\\

Toward this end, it is fairly straightforward to see that a solution to the problem is equivalent to solving the equation $$ax^3+by^3=cz^3$$ in integers for which $rad(a)|x$, $rad(b)|y$, $rad(c)|z$, and $\gcd(ax,by)=1$. We will make things as simple as possible for ourselves by restricting our attention to the case $a=b=1$ and $rad(c)=p$, an odd prime, which results in the problem of solving one of
$$x^3+y^3=p^{\mu}z^3 \; \; \; \; (\mu \in {1,2})$$
with $(x,y)=1$ and $p|z$. We will focus primarily on the case $\mu = 1$, as the other case is quite similar. The following result provides a solution to Erd\"{o}s' problem, however it seems to be merely a first step in this direction.\\

\begin{thm}
Let $p$ denote an odd prime for which the curve 
$$E: Y^2=X^3-432p^2$$
has positive rank. Then there are infinitely many pairwise coprime integer solutions $(x,y,z)$ to $x^3+y^3=p^4z^3$. Furthermore, if $P$ denotes a generator of infinite order on $E$, then pairwise coprime integer solutions to $x^3+y^3=p^4z^3$ can be derived from $(3pk)P$ for every integer $k \ge 1$.
\end{thm}

To be more explicit we elucidate the above statement. Assume that $u,v,d$ are integers for which $(uv,d)=1$ and $(X=u/d^2,Y=v/d^3)$ is a point on the curve. Assume further that $p$ divides $d$. We leave it for the reader to verify that the triple $(x,y,z)$ given by\\
$x=num((36pd^3+v)/6ud)$,\\
$y=num((36pd^3-v)/6ud)$,\\
$z=denom((36pd^3+v)/6ud)$\\
(all in lowest terms) gives solutions to Erd\"{o}s' problem.\\

It is worth remarking at this point that in order to ensure $(x,y)=1$ and $p$ divides $z$, a sufficient condition is for $p$ to divide the integer $d$ defined just above. However, this condition is often not actually necessary. This corresponds to the extra factor of $3$ in front of $P$ in the statement of the theorem, which seems to be necessary only when $p \equiv 1 \; (\bmod \; 3)$.\\

Adam Logan has pointed out that the $3$-Selmer group of $E$ has order $1$ for $p \equiv 4,7,8 \; (\bmod \; 9)$. Consequently, one would expect the rank of $E$ to be $1$ for all of these cases.\\

\section{Proof of Theorem 1.1}

We begin by connecting the $p$-divisibility of $d$ to that of $z$. The definitions of $x$, $y$ and $z$ from the previous section imply the existence of a non-zero integer $k$ for which
$36pd^3+v=kx, 36pd^3-v=ky, 6ud=kz.$ Assuming that $p|d$, then because of the fact that $u$ and $v$ are coprime to $d$, it follows that $p$ does not divide $kx$, and therefore cannot be a factor of $k$. The equation $6ud=kz$ shows that $p$ must be a divisor of $z$.\\

We therefore need to pin down the set of points $(X,Y)=(u/d^2,v/d^3)$ on $E: Y^2=X^3-432p^2$ which have the property that $p$ divides $d$. Let $P$ denote a point of infinite order on $E$, then regarded as a point on $E(\mathbb{Q}_p)$, the multiple of $P$ having denominator divisible by $p$ is equivalent to the order of $P$ in $E(\mathbb{Q}_p)/E_0(\mathbb{Q}_p)$, where $E_0(\mathbb{Q}_p)$ consists of those points whose reduction is non-singular. Since $E$ has additive reduction, this quotient consists of the additive group $\mathbb{Z}/p\mathbb{Z}^{+}$ of order $p$ times the order of the group of components of the Neron model (see p.359 of \cite{S}), which in this case is $3$. It follows that the point $(3p)P$ has coordinates $(u/d^2,v/d^3)$ satisfying the property that $d$ is divisible by $p$, which translates into the corresponding value for $z$ in the solution to $x^3+y^3=pz^3$ being divisible by $p$, and thus giving integer solutions to Erd\"{o}s' problem on $3$-powerful numbers.\\

It is worth noting that if $m$ is squarefree, the additive reduction at each prime can be pieced together so that the group structure coming from each prime dividing $m$ can also be glued together to get pairwise coprime integer solutions to $x^3+y^3=m^4z^3$ by multiplying a generator by $3m$. A simple example of this phenomenon is given by $m=35$, which we leave to the reader for verification.\\

\section{An Elliptic Curve Analogue of the Ankeny-Artin-Chowla Conjecture}

In 1952, Ankeny, Artin and Chowla \cite{AAC} published a congruence involving various quantities related to a quadratic field and a Bernoulli number, which resulted in a statement which is now referred to as the AAC conjecture. Specifically, the AAC conjecture states that if $p \equiv 1 \; (\bmod \; 4)$ is prime, and if $\epsilon_p=\frac{t+u\sqrt{p}}{2}$ is the fundamental unit of the quadratic field $K=\mathbb{Q}(\sqrt{p})$, then $p \nmid u$. Mordell \cite{Mo} later conjectured this to be true for $p\equiv 3 \; (\bmod \; 4)$.\\ 

An equivalent formulation of this can be stated as follows. If $u=\frac{a+b\sqrt{p}}{2} > 1$ is a unit in the ring of integers of $K$ with $b \equiv 0 \; (\bmod \; p)$, then $u=\epsilon^{kp}$ for some integer $k$.\\

This last statement appears to have similarities with the required properties that were required to solve $x^3+y^3=p^4z^3$ in pairwise coprime integers.\\

Let $p$ denote any odd prime, and let $E_{p,2}: Y^2=X^3-432p^2$ and $E_{p,4}: Y^2=X^3-432p^4$.
Let $P=(X,Y)$ denote a generator of one of the curves $E_{p,2}$ or $E_{p,4}$, and $Q=kP=(u/d^2,v/d^3)$ for integers $k,u,v,d$. We have seen that if the denominator $d$ is divisible by $p$, then $Q$ gives rise to a solution in pairwise coprime integers $(x,y,z)$ of
$x^3+y^3=p^4z^3$ and $x^3+y^3=p^5z^3$ respectively. Although Theorem 1.1 states that the integer $k$ should be divisible by $3p$, our computations have shown that the desired condition holds very often when $k$ is divisible by $p$, and thus we state the following as an analogue of AAC for these two families of elliptic curves.\\

\begin{conj} EC-AAC\\
If $P$ is a generator of $E_{p,2}$ (resp. $E_{p,4}$), and $k$ is a positive integer for which $Q=kP=(u/d^2,v/d^3)$ with $p|d$, then $p$ divides $k$.
\end{conj}

A few cautionary historical remarks are in order. At the time that this author learned of the AAC conjecture, the obvious question of whether a similar result holds for composite discriminants led to the immediate discovery that the fundamental units $\frac{t+u\sqrt{d}}{2}$ in $\mathbb{Q}(\sqrt{d})$, with $d=46$, $430$ and $1817$, actually have the property that $d|u$. This led to an investigation by Stephens and Williams \cite{SW}, who produced a longer list of composite values with this property, and very recently, new examples have been found by Reinhart \cite{Re}, with one of the examples being a prime of the form $4k+3$, and hence a counterexample to Mordell's extension of AAC. Therefore, one should be extremely wary about the truth of the AAC conjecture. Similarly, the composite value $m=1349$ has the property that $E_{m,4}$ is generated by a point $P$ having $1349|d$, and presumably a larger search would produce other such curves. We have yet to find a composite integer $m$ for which $m$ divides the integer $d$ arising from a generator of the curve $E_{m,2}$, nor have we found a counterexample to Conjecture 3.1.

\noindent {\bf Acknowledgements} The author would like to express gratitude to Noam elkies, Adam Logan and Joseph Silverman for their extremely helpful discussions.

\end{document}